\documentclass[a4paper,twoside,11pt]{amsart}

\usepackage{amssymb}
\usepackage{amsmath}
\usepackage{epsfig}

\usepackage[latin1]{inputenc} 
\usepackage[T1]{fontenc}

\newtheorem{The}{Theorem}
\newtheorem{Lem}{Lemma}
\newtheorem{Cor}[The]{Corollary}
\newtheorem{Pro}[Lem]{Proposition}
\newtheorem{Hyp}{Hypothesis}

\def\a{\alpha}

\def\C{\mathbf C}
\def\D{\Delta}
\def\d{\delta}

\def\f{\phi}

\def\l{\lambda}

\def\O{\mbox{\rm orb}}
\def\p{\pi}

\def\R{\mathbf R}
\def\r{\rho}

\def\s{\sigma}
\def\t{\tau}

\def\Z{\mathbf Z}
\def\z{\zeta}

\title{Bijectiveness of the Nash map for quasi-ordinary hypersurface singularities}
\author{P.D. Gonz\'alez P\'erez}

\address{Departamento de Algebra. Facultad de Ciencias Matem\'aticas. Universidad Complutense de Madrid.
Plaza de las Ciencias 3. 28040. Madrid. Spain.}

\email{pgonzalez@mat.ucm.es}


\thanks{The authors' research was supported in part by {\em Programa Ram\'on y Cajal} and  MTM2004-08080-C02-01 grants of {\em Ministerio de Educaci\'on y Ciencia}, Spain.}

\begin{document}

 \begin{abstract} In this paper we give a positive answer to a question of Nash, concerning the arc space
 of a singularity,  for the class of quasi-ordinary hypersurface singularities extending to this case previous
  results and techniques of Shihoko Ishii.
\end{abstract}

\maketitle

\section*{Introduction}

In a 1968 preprint,  published later as
\cite{Nash},  Nash introduced a map, nowadays called the {\em Nash map}
(resp. {\em local Nash map}), from
the set of families of arcs with origin at the singular locus of a
variety $X$ (resp. at a fixed singular point $x \in X$),
to the set of essential divisors over the singular
locus of $X$ (resp. over the point $x$).
These families of arcs are called {\em Nash components} (resp. {\em local Nash components}).
Obviously, these maps coincide if the singular locus of $X$ is reduced
to an isolated singular point $x$. Nash showed that these maps
are injective and asked if they were surjective.
Ishii and  Koll\'ar  have shown an affine four
dimensional variety with singular locus reduced to a point, for which the
Nash map is not bijective (see \cite{IK}).
The answer of Nash question for surface and threefolds
singularities is not known in general.  Plenat has given sufficient conditions for the surjectivity of the Nash map for
isolated surface singularities in \cite{Plenat1}.
In the case of surface
singularities the Nash map it is known to be bijective in the following cases: for minimal singularities by Reguera \cite{Reguera},
for  sandwiched
singularities by the work of Lejeune and Reguera \cite{LR} and \cite{Reguera2},
for rational double points of type $A_n$, already studied by Nash \cite{Nash}, and of type $D_n$ by  Plenat \cite{Plenat2} (the result for rational double points in general is announced in  \cite{Plenat}). Plenat and
Popescu-Pampu have shown a class
varieties of dimension two and higher
for which the Nash map is
bijective in \cite{PPlenat} and \cite{PPlenat2};
a similar result in the surface case is
announced by Morales \cite{Morales}.

Ishii and Koll\'ar have shown that Nash question
has a positive answer in the case of normal toric
varieties, see \cite{IK}.
 Ishii has generalized this result for the class
pretoric algebraic varieties, which contains in particular the
class of toric varieties (non necessarily normal). Petrov
formulated {\em Nash question for pairs $(X, B)$}, consisting of a
variety $X$ and a proper closed subvariety $B$ containing the
singular locus of $X$, and exhibited a positive answer in the case
of pairs, $(X, B)$, formed by  a normal toric variety $X$ and an
invariant set $B$. He applied this result to prove the
bijectiveness of the Nash map for the class of {\em stable toric
varieties}, a class of reduced but non necessarily irreducible
varieties introduced by Alexeev \cite{Alexeev}, which generalize
normal toric varieties (see \cite{Petrov}).
 Ishii has
shown recently that the local Nash map is bijective for {\em
analytically pretoric} singularities, a class of singularities
containing toric and analytically irreducible quasi-ordinary
hypersurface singularities (see \cite{Ishii-fourier}). The class
of quasi-ordinary singularities appears classically in Jung's
strategy to obtain resolution of surface singularities from the
embedded resolution of plane curves (see \cite{Jung},
\cite{Abhyankar} and \cite{Lipman2}).

The purpose of this Note is to show that the Nash map is bijective
for any reduced germ of quasi-ordinary hypersuperface singularity,
Theorem \ref{main}. We show that the results and approach for the
class of {\em pretoric singularities}, analyzed by Ishii in
\cite{Ishii-crelle}, can be extended to the case of an
analytically irreducible quasi-ordinary hypersurface germ $(X,x)$
by showing a property, Proposition \ref{llave}, deduced from
Lipman's result on the structure of the singular locus (see
\cite{Lipman}, \S 7). Moreover, if $B$ is a proper closed subgerm
of $(X,x)$ containing the singular locus of $X$, we prove that
Nash question for the pair $(X, B)$ has a positive answer under a
natural technical condition (see Hypothesis \ref{Hyp} and Theorem
\ref{na}).

In the case of a quasi-ordinary hypersurface germ  $(X, x)$ with
several irreducible components $X_i$, $i =1, \dots, t$, the
analysis of the pre-image of the singular locus of $X$ by the
normalization map (see Lemma \ref{nu}) is essential to deduce the
main result Theorem \ref{main}, the bijectiveness of the Nash map
for $X$, from the bijectiveness of the Nash map for suitable pairs
$(X_i, B_i)$ for  $i =1, \dots, t$.

Finally, we compare the notions of pretoric singularity and
analytically pretoric singularity, introduced by Ishii in
\cite{Ishii-crelle} and \cite{Ishii-fourier} respectively, with
                           the notion of {\em toric quasi-ordinary singularity} introduced
\cite{GP03} and we show that our main result holds in a slightly larger category, which we call
{\em strongly analytically pretoric} (see Corollary \ref{na-extension}).

The explicit description of the essential divisors over the singular locus of a
quasi-ordinary hypersurface singularity,  given in this paper, is
applied by Hernando and the author in \cite{GP-H} to prove that the essential
divisors of an irreducible germ $(X, x)$  of  quasi-ordinary hypersurface
determine, through a suitable notion of Poincar\'e series, the
{\em characteristic monomials}, an analytical invariant which in the
analytic case encodes the embedded topological type of the germ
characterized by the work of Gau and Lipman (see \cite{Lipman} and \cite{Gau}).
It should be noticed that the Poincar\'e series associated only to the
essential divisors over the point $x$, which correspond by Ishii's result
\cite{Ishii-fourier} to the local Nash components, does not contain enough information to recover
the characteristic monomials of the quasi-ordinary hypersurface $(X, x)$ in general, see \cite{GP-H}.  This property reflects the fact that
quasi-ordinary singularities are rarely isolated: $X$ does have an isolated singularity at $x$ if and only if the singular germ $(X,x)$  is of dimension one or normal of dimension two.
The use of all essential divisors, over the different components of the singular locus of $X$, and also over $x$,  is crucial to recover
the characteristic monomials of the germ $(X,x)$, as main local invariants.

\section{Basic definitions on the arc space and relative Nash
problem} \label{pairs}

In this section we give basic definitions and results on the {\em
relative Nash problem}, also called {\em Nash problem for a pair}.
These notions are natural extensions of the corresponding ones on
the classical Nash problem (see \cite{Petrov}). In this paper the
scheme $X$  is a pure dimensional reduced algebraic (resp.
algebroid germ of)  variety, defined over  a field $k$,
algebraically closed of characteristic zero. Let $B \subset X$ be
a reduced proper $k$-subscheme containing the singular locus
$\mbox{\rm Sing} (X)  $ of $X$.

 A {\em resolution of singularities of $(X, B)$} is a proper
modification $\f: Y \rightarrow X$ such that $\f_{|Y -\f^{-1} (B)
}: Y - \f^{-1} (B) \rightarrow X - B$ is an isomorphism. The
resolution $\f$ is divisorial is $\f^{-1} (B)$ is a divisor. An
{\em exceptional divisor $E$ over $X$ relative to $B$} is an
exceptional divisor such that the center of $E$ over $X$ is
contained in $B$. An exceptional divisor $E$ over $X$ relative to
$B$ is {\em essential} if for every resolution $\f : Y \rightarrow
X$ of the pair $(X,B)$ the center of $E$ on $Y$ is an irreducible
component of $\f^{-1} (B)$. This center is called an {\em
essential component} on $Y$.

If $k \subset K$ is a field extension an {\em arc} over $X$ is a
morphism $\a: \mbox{\rm Spec} K[[t]] \rightarrow X$. We denote
respectively  by $0$ and $\eta$  the closed point and generic
point of $\mbox{\rm Spec} K[[t]]$. If $m \geq 0$ is an integer, an
$m$-{\em jet} over $X$ is a morphism $\a: \mbox{\rm Spec}
K[[t]]/(t^{m+1}) \rightarrow X$.  The set $X_m$ of $m$-jets can be
given the structure of scheme of finite type over $k$. We have
canonical morphisms  $X_{m+k} \rightarrow X_{m}$, corresponding to
$\mbox{\rm Spec} K[[t]]/(t^{m+1})  \rightarrow \mbox{\rm Spec}
K[[t]]/(t^{m+k+1})$,   for all $m , k \geq 0$.
 The
arc space $X_{\infty} := \underleftarrow{\lim} X_m $ has  the
structure of scheme over $k$, not of finite type.  A point $z \in
X_{\infty}$ corresponds to an arc $\a_z :  \mbox{\rm Spec} K[[t]]
\rightarrow X$ such that $K$ is the residue field at $z$. We have
that $K$-valued points of $X_{\infty}$ correspond to arcs of the
form $\mbox{\rm Spec} K[[t]] \rightarrow X$ bijectively. We often
denote the point $z$ and the corresponding arc $\a_z$ by the same
symbol. We have a canonical projection $ \p^X: X_{\infty}
\longrightarrow X$, defined by $ z \mapsto \a_z (0)$. A morphism
of varieties (resp. of algebroid germs) $\psi: Y \rightarrow
X$, corresponds to a morphism at the level of arcs: $
\psi_{\infty} : Y_{\infty} \rightarrow X_{\infty}$. See \cite{DL} for
the general construction of arcs spaces.

A {\em Nash component of $X$ relative to $B$} is an irreducible
component of $\p_X^{-1} (B)$ which is not contained in $B_\infty
\subset X_{\infty}$. If the field $k$ is of characteristic zero, as it is assumed in this paper,
then the proof of Lemma 2.12 \cite{IK} imply that the Nash components of $X$ relative to $B$ coincide with
the irreducible components of $\p_X^{-1} (B)$.
 Denote by $\bigcup_{i} C_i$ the union of Nash
components of $X$ relative to $B$. Let $\f: Y \rightarrow X$ be a
divisorial resolution of $(X, B)$, i..e, $\f^{-1} (B)$ is a
divisor with irreducible components $E_1, \dots, E_l$. Then the
restriction $\f_{\infty}$ to $ \bigcup_{j=1}^l \p_Y^{-1} (E_j)
\longrightarrow \bigcup_{i} C_i$. is dominant and bijective
outside $B_{\infty}$. For all $i$ there exists a unique $j_i$ such
that $\p_Y^{-1} (E_j) \rightarrow  C_i$ is dominant. The analogous
statement of Nash Theorem in this relative situation is that
$E_{j_i}$ is an essential divisor of $X$ relative to $B$  and that
the {\em relative Nash map} $ C_i \mapsto E_{j_i} $ is an
injection, from the Nash components and to the set of essential
divisors of $X$ relative to $B$ (see Lemma 2.14
\cite{Ishii-fourier} for an analogous proof in this relative
situation or Theorem 2.17 \cite{Petrov} for a sketch of proof in
the algebraic case). See also \cite{IK} for a modern proof of the
classical statement of Nash \cite{Nash}, when $B = \mbox{Sing}
(X)$). The Nash problem for $X$ relative to $B$, i.e., the Nash
problem for the pair $(X, B)$,  is to determine if this
correspondence is bijective.

\section{Basic notations on normal toric varieties}

We give some basic definitions and notations (see \cite{Oda}
or  \cite{Fulton} for proofs). If $N$ is a lattice we denote by
$M$ the dual lattice, by $N_\R$ the real vector space spanned by
$N$ and by $\langle , \rangle$ the canonical pairing between the
dual lattices $N$ and $M$ (resp. vector spaces $N_\R$ and $M_\R$).
A {\it rational convex polyhedral cone} $\tau $ in $N_\R$, a {\it
  cone} in what follows, is the set $\t := \mbox{\rm pos} \{a_1,
\dots, a_s\}$
of non negative linear
combinations of vectors $a_1, \dots, a_s \in N$. The cone $\tau$
is {\it strictly convex} if $\tau $ contains no linear subspace
of dimension $>0$. We denote by $\stackrel{\circ}{\tau}$ the {
\it relative interior} of a cone $\sigma$. The {\it dual} cone
$\tau^\vee$
(resp. {\it orthogonal} cone $\tau^\bot$)
of
$\tau$ is the set $ \{ w  \in M_\R / \langle w, u \rangle \geq
0\}$ (resp. $ \langle w, u \rangle = 0$)  $ \; \forall u \in \tau \}$).
A {\it fan\index{fan}} $\Sigma$ is a family of
strictly convex
  cones  in $N_\R$
such that any face of such a cone is in the family and the
intersection of any two of them is a face of each.
If  $\tau$ is a cone in the fan $\Sigma$, the
semigroup $\tau^\vee  \cap M$ is of finite type, it spans the
lattice $M$ and the variety
$Z_{\t, N} = \makebox{Spec} \, k [ \tau^\vee \cap M ]$,
which we denote by $Z_\tau$ when the lattice
is clear from the context, is normal.
The affine varieties $Z_\tau$
corresponding to cones in a fan $\Sigma$ glue up to define the
{\it toric variety\index{toric variety}} $ Z_\Sigma$. The torus $ T_N
:= Z_{ \{ 0 \} } \cong
(k^*)^{ \mbox{rk}\, N}  $ is embedded in $ Z_\Sigma$ as an open dense subset
and there is an action of $T_N$ on $Z_\Sigma$ which extends the action
of
the torus on itself by multiplication.  We have a bijection between the
relative interiors of the cones  of the fan and the orbits of the
torus action, $\stackrel{\circ}{\tau} \mapsto \O_{Z_\Sigma} (\tau)  $,
which inverses inclusions of the closures. We denote the orbit
$\O_{Z_\Sigma} (\tau)$ by $\O (\tau)$ when the toric variety
$Z_\Sigma$ is clear from the context.

In this paper $\s$ denotes a rational strictly convex cone in $N_\R$ of maximal dimension.
Any non zero vector $v \in \s \cap N$ defines a
valuation $  \mbox{\rm val}_{v}  $ of the field of fractions of $k[\s^\vee \cap M ]$ (resp. of the $m$-adic completion  $k[[\s^\vee \cap M ]]$ of the localization of $k[\s^\vee \cap M ]$ at the maximal ideal, $(\s^\vee \cap M) \setminus \{ 0\}$, defining the {\em origin} $o_\s$
of the toric variety $Z_{\s}$). This valuation, called {\em monomial} or {\em toric valuation}, is
defined for an element $0 \ne \f = \sum c_u X^u \in k [[ \s^\vee \cap M
]]$, by $ \mbox{\rm val}_{v}  (\f) = \min_{c_u \ne 0} \langle n, u \rangle$.
If the ray $\r:= v \R_{\geq 0}$ belongs a fan $\Sigma$ subdividing $\s$,
the closure $D_\r $ of the orbit  ${\O (\r) }$ is an invariant divisor
(we denote it also by $D_v$ if the vector $v$ is primitive for the
lattice $N$).
We denote by  $\mbox{\rm val}_{D_{\r}}$ the associated divisorial valuation.
If $v= q v_0$ for $q \in \Z_{\geq 1}$ and   $v_0$ a primitive vector
for the lattice $N$ then we have that  $ \mbox{\rm val}_{v}
=   q \mbox{\rm val}_{D_{\r}}$.  Following Ishii we say that the valuation $\mbox{\rm val}_{v} $ is a {\em toric divisorial valuation} (see \cite{Ishii-crelle}). The cone $\s$ induces a partial order on $N_\R$, defined by $u \leq_\s u'$ iff $u'\in u + \s$ (similar definition for $\leq _{\s^\vee}$ on $M_\R$ holds).

\section{Quasi-ordinary hypersurfaces: singular locus and normalization}

 A quasi-ordinary hypersurface singularity  $X$ is defined by $\mbox{\rm
Spec} \, k [[ {x} ]][y]/(f)$ where  ${x}= (x_1,
\dots, x_d) $ and
$f \in k [[ {x} ]][y]$ is a  {\em quasi-ordinary polynomial}, i.e., a
Weierstrass polynomial such that the discriminant $\D_y f$ with
respect to $y$ is of the form $x^{\d} \epsilon$ where $\d \in
\Z^{d}_{\geq 0}$ and $\epsilon \in k[[ {x}
]]$ is a unit. By the Jung-Abhyankar Theorem the
roots of quasi-ordinary polynomials are fractional power series
with particular properties.
Namely, if $f$ is irreducible of degree $n$ a root of $f$ is of
the form:
$
\z = \sum c_\l x^{\l} \in k [[ {x}^{1/n}]]$,
where ${x}^{1/n} = (x_1^{1/n}, \dots, x_d^{1/n})$, and
the terms appearing in this expansion verify certain properties. In particular
if $f$ is  of degree $> 1$ in $y$,  certain monomials, determined
by comparing the different roots of $f$ and  called {\em characteristic or distinguished},  appear in the
expansion of $\z$ with non zero coefficient. The corresponding exponents, which are also called characteristic,
can be relabelled in the form $ \l_1 \leq_{\s^\vee} \l_2
\leq_{\s^\vee} \cdots \leq_{\s^\vee} \l_g $  . These
exponents determine the following nested sequence of lattices: $
M_0 \subset M_1 \subset \cdots \subset M_g =: M $ where $M_0
:=\Z^{d} $ and $M_j := M_{j-1} + \Z \l_j$ for $j=1, \dots, g$ with
the convention $\l_{g +1} = + \infty$.
We have that the exponents appearing in the expansion of $\z$
belong to $M$. See \cite{Lipman} and also \cite{GP03}.
We have ring extensions:
\begin{equation} \label{extension}
k [[ \s^\vee \cap M_0 ]] = k [[{x}]]
\longrightarrow \mathcal{O}_X \cong k [[ {x} ]][\z]
\longrightarrow k [[ \s^\vee \cap M]]
\end{equation}
where $\s^\vee$ denotes the positive quadrant $\R^{d}_{\geq 0}$
and $M_0 = \Z^d$ (we denote by $\s$, $N_0$ and $N$ the
corresponding dual objects of $\s^\vee $, $M_0$ and $M$ respectively). In \cite{GP03}
it is proved that the ring extension $\mathcal{O}_X \rightarrow k
[[ \s^\vee \cap M]]$ is the inclusion of $\mathcal{O}_X$ in its
integral closure in its field of fractions. Geometrically,
(\ref{extension}) corresponds to a sequence of finite maps:
\begin{equation} \label{integral}
(X_{\s, N} ,o_\s) =(\bar{X},x) \stackrel{\nu}{\longrightarrow}
(X,x) \stackrel{\r}{\longrightarrow} (Z,x) = (X_{\s, N_0} ,o_\s) =
(k^d, 0).
\end{equation}
Since the
map $\r \circ \nu$ is equivariant, it maps the orbit $\mbox{\rm
orb}_{\bar{X}} \t$ to $\mbox{\rm orb}_{{Z}} \t$, for each face $\t
< \s$.

We recall Lipman's description of the singular locus of a
quasi-ordinary hypersurface, see Theorem7.3  \cite{Lipman},
for a precise statement (cf. also the reformulation of this result given in \cite{PP}).

\begin{The}  \label{lip}
With the previous notations if $(X, x)$ is analytically
irreducible we have that the irreducible components of $\mbox{\rm
Sing} (X)$ are of codimension one or two. The codimension one
(resp. two)  components are intersections of $X$ with  $x_i =0$,
(resp. with $x_i = 0$ and $x_j =0$)  for some suitable coordinate
sections, in each case,  determined by the characteristic monomials. $\Box$
\end{The}

\begin{Pro} \label{llave}
The set $\nu^{-1} (\mbox{\rm Sing} (X) ) $ is a germ of
closed set at the origin of $ \bar{X}$, which is invariant by the torus action on $\bar{X}$.
\end{Pro}
{\em Proof.} By Lipman's theorem, the irreducible components of $\mbox{\rm
Sing} (X)$ are the germs   $\r^{-1} ( {\overline{ \mbox{\rm
orb}_{{Z}} (\t)}} )$ at the point $x$, where $\t$ runs certain set
of a one (resp. two ) dimensional faces of $\s$. It follows from
this and the previous discussion that the irreducible components
of $\nu^{-1} ( \mbox{\rm Sing} (X) ) $ are of codimension one or
two. The codimension one (resp. two) components are of the form
\[ (\r \circ \nu)^{-1} (
{\overline{ \mbox{\rm orb}_{{Z}} (\t)}} ) = {\overline{ \mbox{\rm
      orb}_{\bar{X}} (\t) }}  .\]
If $x_i = X^{u_i} $ for $i=1, \dots, d$,   in (\ref{extension}) then
$x_i =0$ (resp. $x_i =x_j =0$) defines in $Z$ the closure of the orbit
$\mbox{\rm orb}_Z (\t)$ where the cone $\t$ is characterized by
$\t^\bot \cap \s^\vee =\mbox{\rm pos}_{k\ne i} ( u_k ) $ (resp. by  $\t^\bot \cap
\s^\vee = \mbox{\rm pos}_{k \ne i, j} ( u_k) $).
$\Box$

\section{Relative Nash problem, the irreducible case} \label{bna}

We follow Ishii's  approach in  \cite{Ishii-crelle} and
\cite{Ishii-fourier}. Let $(X,x)$ be an irreducible germ of
quasi-ordinary hypersurface. We study the relative Nash problem
for a proper closed set $ B \subset X$ such that $\mbox{\rm Sing} (X)
\subset B$. We introduce the following hypothesis on $B$.

\begin{Hyp} \label{Hyp}
 Any irreducible component of $ \nu^{-1} ( B )
$ is an orbit closure ${\overline{ \mbox{\rm orb}_{\bar{X}} (\t)
}} $ corresponding to some face $\t$ of $\s$.
\end{Hyp}
\medskip

\noindent Notice that $\mbox{\rm Sing} (X)$ verifies this
hypothesis by Proposition \ref{llave}. When the hypothesis above
is verified any irreducible component in the closure of the set
$
  \nu^{-1} (B) - \mbox{\rm Sing} (\bar{X})
$
is an orbit closure corresponding to some regular face $\t < \s$
(with respect to the lattice $N$). We denote by $\t_1, \dots,
\t_r$ the regular faces determined in this way and by $e_i \in N$,
the barycenter of $\t_i$ (i.e., the sum of the primitive integral
vectors, for the lattice $N$, in the edges of $\t_i$), for  $i=1,
\dots,r$.

 Let
$\{v_j \}_{j=1}^s $ the set of minimal elements, with respect of
the partial order $\leq_\s$ in the set:
\[
S := \bigcup_{\t < \s, \, \t \mbox{ singular}} \stackrel{\circ}{\t}
\cap N.
\]

By \cite{IK} the toric divisors $\{ D_{v_j}\}_{j=1, \dots,
s} $ are the essential divisors of toric variety $\bar{X}$, and also
the essential divisors of the germ of $\bar{X}$ at
the closed orbit, by Lemma 4.9
\cite{Ishii-fourier}. This characterization of essential divisors
generalizes a result of Bouvier \cite{Bouvier}, see also  \cite{B-GS}.
\begin{Lem} \label{min}
Each $e_i $, for $i= 1, \dots, r$ is minimal among $\{ e_i, v_j
\}_{i=1, \dots, r}^{j=1, \dots, s}$ for the order $\leq_\s$.
\end{Lem}
{\em Proof.} See  Lemma 5.7 \cite{Ishii-crelle}. $\Box$

\begin{Lem} \label{comp}
Let $\{ e_i, v_j \}_{i=1, \dots, r}^{j=1, \dots, w}$, ($w \leq s$)
be the set of minimal elements of  $\{ e_i, v_j \}_{i=1, \dots,
r}^{j=1, \dots, s}$. Then there is an inclusion
\[
\{ \mbox{essential divisors over } X \mbox{ relative to } B \}
\subset \{ D_{e_i}, D_{v_j}\} _{i=1, \dots, r}^{j=1, \dots, w}.
\]
\end{Lem}
{\em Proof.} The statement can be translated in purely
combinatorial terms in terms of the existence of resolutions of singularities of
$X$, which are obtained by composing the normalization map with
toric modifications. The precise arguments are the content of the
proof of Lemma 5.7 \cite{Ishii-crelle}. $\Box$

Following Ishii, \cite{Ishii-crelle},  we associate to a non zero
vector  $v \in \s \cap N$  a subset $T^X_{\infty} (v)$ of  $X_{\infty}$ , containing only arcs
which lift to  arcs
with generic point in the torus $T_N $ of $\bar{X}$:
\[
T^X_{\infty} (v) := \{ \a \in X_{\infty} \mid \a(\eta) \in \nu (T_N ),
\mbox{\rm ord}_t \a^* (x^u) = \langle v, u \rangle, \, \mbox{ for } u \in
M \}.
\]
The sets $T^{\bar{X}}_{\infty} (v)$, defined similarly, are orbits
of a natural action of $(T_N)_\infty$ on the arc space of the normal
toric variety $\bar{X}$ (see \cite{Ishii-algebra}).

\begin{Lem} \label{comp2}
Let $\{ e_i, v_j \}_{i=1, \dots, r}^{j=1, \dots, w}$, ($w \leq s$)
be the set of minimal elements of  $\{ e_i, v_j \}_{i=1, \dots,
r}^{j=1, \dots, s}$. Then, the following closures are distinct Nash
components of $X_\infty$:
\[
\overline{T^X_{\infty} (e_i) },\, i =1, \dots, r\,  \mbox{ and }
\, \overline{T^X_{\infty} (v_j) }, j=1, \dots, w.
\]
If $v \in \{ e_i, v_j \}_{i=1, \dots,
  r}^{j=1, \dots, s}$ the image of the component
$\overline{T^X_{\infty} (v) }$
by the Nash map is the divisor $D_v$.
\end{Lem}
{\em Proof}.  The proof is analogous to the proofs of Lemma 4.6
and 4.7 in \cite{Ishii-fourier}. Notice that  with our hypothesis
the proof holds not only for vectors $v \in \stackrel{\circ}{\s}
\cap N $ but  also for vectors $v \in {\s} \cap N $. See also the
proof of Lemma 5.10 and Theorem 5.11 \cite{Ishii-crelle}.
$\Box$

\begin{The} \label{na}
Let $(X,x) $ be a irreducible germ of quasi-ordinary hypersurface
singularity.  Let $B$ a closed subscheme verifying Hypothesis \ref{Hyp}.
Then the Nash map between the set of Nash components of
$\p_X^{-1} (B) $ and the set of essential divisors of $X$ relative
to $B$ is bijective.
\end{The}
{\em Proof}. Let  $u \in \{ e_i, v_j \}_{i=1, \dots, r}^{j=1,
\dots, w}$. Then the sequence
\[
u \mapsto \overline{T^X_{\infty} (u) } \mapsto D_u
\]
defines an injection from the set $\{ e_i, v_j \}_{i=1, \dots,
r}^{j=1, \dots, w}$  to the set of essential divisors over $X$, by
Lemma \ref{comp2} and by the injectivity of the Nash map. By Lemma
\ref{comp} the set of essential divisors is of cardinality less or
equal than $r+w$, hence it follows that this set is of cardinality
equal to $r+w$ and the injection above is a bijection. This
implies that the Nash map is bijective. $\Box$

\begin{Cor}
If $(X,x)$ is an analytically irreducible
quasi-ordinary hypersurface  the Nash map is bijective.

\end{Cor}
{\em Proof}.  By Proposition \ref{llave} the singular locus $B
= \mbox{\rm Sing} (X)$ verifies the hypothesis \ref{Hyp} of Theorem \ref{na}.$\Box$

\section{Nash problem for a quasi-ordinary hypersurface} \label{na2}

Now we suppose that $(X, x)$ is a germ of reduced quasi-ordinary
hypersurface. We denote by $f$ a {\em quasi-ordinary polynomial}
defining $(X,x)$. The factors $f_i$ of the factorization of $f=
f_1\dots f_t$ as product of irreducible terms corresponds to the
irreducible components of $(X,x)$. The factors $f_i$ are
quasi-ordinary polynomials, for $i=1, \dots, t$.

We denote by $B_i $ the intersection:
\[
B_i = X_i \cap \mbox{\rm Sing} (X) = \mbox{\rm Sing}(X_i) \cup
\bigcup_{j=1, \dots, t} ^{j\ne i} X_i \cap X_j.
\]

We denote by $\nu_i : \bar{X}_i \rightarrow X$ the normalization
of $X_i$, which is a toric singularity by the previous discussion.

\begin{Lem} \label{nu}
We have that $\nu_i^{-1} (B_i) $  is a germ of invariantly closed
set,  at the close orbit of the toric singularity $\bar{X}_i$.
\end{Lem}
{\em Proof}. We have already shown the statement for $\nu_i^{-1}
(\mbox{\rm Sing}(X_i))$ by Proposition \ref{llave}. If $j \ne i$
then  $\nu^{-1} (X_i \cap X_j)$ is defined by $f_j (\z^{(i)})=0$
where $\z^{(i)}$ is a fixed root of $f_i$, used to define  the
ring extension (\ref{integral}) corresponding to $X_i$.
We have that the element $f_j (\z^{(i)})$ is equal to the product
of  a monomial by a unit in the local algebra of the toric
singularity $\bar{X}_i$ (this follows easily from the definition,
see \cite{GP03} for more details). This implies that $f_j
(\z^{(i)}) =0$ defines
 a germ of invariantly closed
set, at the close orbit of the toric singularity $\bar{X}_i$,
which is equal to $\nu_i ^{-1} (X_i \cap X_j)$. $\Box$

We obtain then a generalization of  Corollary 4.12 in
\cite{Ishii-fourier}:

\begin{The} \label{main} Let $(X,x)$ be a reduced  germ of quasi-ordinary hypersurface
singularity. Then  the Nash map between the Nash components of
$\p_X ^{-1} ( \mbox{\rm Sing} (X)) $ and the essential divisors of
$X$ is bijective.
\end{The}
{\em Proof.} We keep the notations given above for the irreducible
components of $X$. Notice that $\p_X^{-1} ( \mbox{\rm Sing} (X) )
= \bigsqcup_{i=1}^t \p_{X_i}^{-1} (B_i)$ by definition of $B_i$.
It follows from this that:
\[
\{ \mbox{Nash components of } X \} \subset \bigsqcup_{i=1}^t \{
\mbox{Nash components of } X_i \mbox{ relative to } B_i \}.
\]
(See the proof of Lemma 4.11
\cite{Ishii-fourier}).
We prove that:
\[
\{ \mbox{essential divisors over } X \} \subset \bigsqcup_{i=1}^t
\{ \mbox{essential divisors over } X_i \mbox{ relative to } B_i
\}.
\]
Let $\f_i : Y_i \rightarrow X_i$ be a resolution of $(X_i, B_i)$.
Then the composite $\f :Y \rightarrow X$ defined by:
\[
Y:= \bigsqcup_{i=1}^t Y_i \stackrel{ \bigsqcup \f_i }{\rightarrow}
\bigsqcup_{i=1}^t X_i \rightarrow X
\]
is a resolution of the pair $(X, \mbox{\rm Sing} (X))$ by the  definition of
$B_i$. Let $E$ be an essential divisor of $X$. The center of $E$
in $Y$ is an irreducible component of
\[
\f^{-1} ( \mbox{\rm Sing} (X)) = \bigsqcup_{i=1}^t \f_i^{-1}
(B_i),
\]
thus an irreducible component of $\f_i^{-1} (B_i) $ for some $i$.
This implies the assertion.

The hypothesis \ref{Hyp} is verified by $B_i$ with respect to
$X_i$ for $i=1, \dots,t$ by Lemma \ref{nu}. By Theorem \ref{na}
applied to the pair $(B_i, X_i)$ the Nash map between the set of
Nash components of $\p_{X_i}^{-1} (B_i) $ and the set of essential
divisors of $X_i$ relative to $B_i$ is bijective.  $\Box$

\section{An extension of the results to a larger category}

In Definition 4.1 \cite{Ishii-fourier} the notion of
analytically pretoric singularity is introduced in the algebroid
category. A germ $(X,x)$ is called in \cite{Ishii-fourier} {\em
analytically pretoric} if there exists a sequence of injective
local homomorphisms:
\[
k[[\s^\vee \cap M_0]] \stackrel{\r^*}{\rightarrow}
\mathcal{O}_{X,x} \stackrel{\nu^*}{\rightarrow} k[[\s^\vee \cap
M]]
\]
such that \begin{enumerate} \item[i.\label{1}]
 $\nu^* \circ \r^* : k [[\s^\vee \cap M_0]]
\rightarrow k[[ \s^\vee \cap M ]] $ is the canonical injection
corresponding to a finite index lattice inclusion $M_0 \subset M$,
\item[ii.\label{2}]   the morphism $\nu: \mbox{\rm Spec} k [[\s^\vee \cap M]]
\rightarrow X$ corresponding to $\nu^*$ is the normalization map,
\item[iii.\label{3}]  the restriction  of $\nu$ to the torus $\mbox{\rm Spec} k
[[\s^\vee \cap M]][M]$ is an isomorphism onto
 its image.
\end{enumerate}
Notice that Ishii introduced the notion of {\em pretoric variety}
in the algebraic category in \cite{Ishii-crelle} Definition 5.1.
The first two conditions for a variety $X$ to be pretoric are the
formulations of axioms i. and ii. above in the algebraic category,
while the third condition above is to be replaced by
\begin{enumerate}
\item[iii'.\label{3prima}]  The closed subset $\nu^{-1} (\mbox{\rm Sing} (X))$ is invariant for the torus action of
$\bar{X}$.
\end{enumerate}
We say that a germ of algebroid singularity is {\em strongly analytically pretoric} if it verifies
conditions   i., ii.,  and iii'.
Notice that in the algebroid case condition iii'.  implies condition
iii in Ishii's definition of analytically pretoric singularity.

\begin{Cor} \label{na-extension}
  If $(X,x)$ is a germ of strongly analytically pretoric singularity then
the associated Nash map is bijective.
\end{Cor}
{\em Proof}. The analysis and results done in the sections \ref{bna} and \ref{na2}
 extends formally for any algebroid singularity  $(X, x)$
  verifying conditions  i., ii.,  and iii'. $\Box$

We end this section by formulating some natural questions which come out from the comparison
of the notions of (strongly) analytically pretoric with that of {\em toric quasi-ordinary
singularity},  introduced in \cite{GP03}  as a suitable generalization of the
notion of quasi-ordinary singularity (in \cite{GP00} the
this class of singularities was restricted to the case of relative hypersurface
germ in $Z _{\s} \times \C $).  A germ of complex analytic
variety $(X,x)$ of pure dimension $d$ is {\em a toric
quasi-ordinary singularity} if there exists:
\begin{enumerate}
\item[a.] An affine  normal toric variety $X_{\s, N_0} =
\mbox{\rm Spec } k[\s^\vee \cap M_0]$, defined by a $d$
dimensional rational strictly convex cone $\s$ for the lattice
$N_0$ (dual of $\s^\vee$ and $M_0$ respectively). \item[b.] A
finite morphism of germs $(X, x) \stackrel{\r}{\rightarrow}
(X_{\s, N_0}, o_\s)$, where $o_\s$ is the closed orbit, which is
unramified over the torus of $X_{\s, N_0}$.
\end{enumerate}
Condition  b. means that for each representative of the
morphism $\r$ there exists an open neighborhood of $o_\s$ such
that $\r$ is unramified over its intersection with the torus.
If the cone $\s$ is simplicial then $(X,x)$ is a quasi-ordinary
singularity:  we reduce to this case by replacing the lattice
$M_0$ by a sublattice $M_0'$ of finite index such that $\s^\vee$
is regular for $M_0'$ and hence $Z _{\s,N_0'} = \C^d$, for $N_0'$
the dual lattice of $M_0'$. Then the normalization $(\bar{X}, x)$
is a germ of toric variety $X_{\s, N} = \mbox{\rm Spec }
\C [\s^\vee \cap M]$ for some lattice $M \supseteq M_0$,  at is
closed orbit, and the composite $\r \circ \nu$ is a germ of toric
equivariant map defined by the change of lattices (see Theorem 5.1
in \cite{PP}, or see \cite{GP03} in the hypersurface case).
The definition of analytically pretoric singularity
can be easily adapted to the complex analytic category.
It is
immediate that if $(X,x)$ is a germ of analytically pretoric singularity then
it is toric quasi-ordinary singularity, since the map $\r$ is then
unramified over the torus by axiom i and iii. It follows from the discussion above that if
the cone $\s$ is simplicial the notions of toric quasi-ordinary
singularity and that of analytically pretoric singularity coincide
(at least in the complex analytic category). It seems quite reasonable that
both notions coincide also in the algebroid category and without
any assumption on simpliciality on the cone $\s$ appearing in both
Definitions.
Is to
be conjectured that for a equidimensional germ of algebroid singularity conditions
i., ii., iii. and i., ii., iii'. are equivalent.

\medskip

\noindent
{\bf Acknowledgement.} The author is grateful to Profesor Shihoko Ishii for stimulating conversations
during the Hayashibara Forum on Singularities 2006 at IH\'ES.  During the preparation of this work the
author has benefited from the hospitality of {\em Institut de Math\'ematiques de Jussieu} in Paris.

{\small
{\sc \quad Pedro Daniel Gonz\'alez P\'erez}
}


\begin{thebibliography}{GGGG}



\bibitem[Al]{Alexeev}
{\sc Alexeev, V.} Complete moduli in the presence of semiabelian
group action.  {\em Ann. of Math.} (2)  {\bf 155}  (2002),  no. 3,
611--708.

{\sc Abhyankar, S.S.}, On the ramification of algebraic functions.
{\em Amer. J. Math.},  {\bf 77}. (1955), 575-592.


\bibitem[A]{Abhyankar}{\sc Abhyankar, S.S.}, On the ramification of
algebraic functions.
{\em Amer. J. Math.},  {\bf 77}. (1955), 575-592.


\bibitem[Bo]{Bouvier}{\sc Bouvier C.}, Diviseurs essentiels, composantes essentielles des vari\'et\'es toriques
singul\`eres, {\em Duke Math. J.} Volume {\bf 91}, No 3, (1998),
609-620.

\bibitem[Bo-GS]{B-GS}{\sc Bouvier C., Gonzalez-Sprinberg G.,} Syst\`eme g\'en\'erateur minimal, diviseurs essentiels et G-d\'esingularisations de vari\'et\'es toriques, {\em T\^ohoku Math. J. }, Volume {\bf 47} (1995),  125-149.


\bibitem[D-L]{DL}
{\sc Denef, J.; Loeser, F.} Germs of arcs on singular algebraic
varieties and motivic integration. {\em Invent. Math.} {\bf
135} (1999), no. 1, 201--232.


\bibitem[Fu]{Fulton} {\sc Fulton, W.}; {\em Introduction to Toric Varieties}; Annals of Math. Studies (131),
Princeton University Press, 1993.
\bibitem[Gau]{Gau}
{\sc Gau, Y-N.}, {\em Embedded Topological classification of
  quasi-ordinary singularities}, Memoirs of the American Mathematical
Society {\bf 388}, 1988.

 \bibitem[GP1]{GP00}{\sc Gonz\'alez P\'erez P.D.}, Singularit\'es
 quasi-ordinaires toriques et poly\`edre de Newton du discriminant,
{\em Canadian J. Math.} {\bf 52} (2), 2000, 348-368.
00
\bibitem[GP2]{GP03}
{\sc Gonz\'alez P\'erez, P.D.}, Toric embedded resolutions of
quasi-ordinary hypersurface singularities, {\em Ann. Inst. Fourier
  (Grenoble)}, 53 (6), (2003), 1819-1881.

\bibitem[GP-H]{GP-H}
{\sc Gonz\'alez P\'erez, P.D., Hernando F.}, Quasi-ordinary
singularities, essential divisors and Poincar\'e series, arXiv:
0705.0603 [math.AG].


\bibitem[I-K]{IK}
{\sc Ishii, S. \& Kollar, J.}, The Nash problem on arc families of
singularities, {\em Duke Math. J.}, {\bf 120}  (2003),  no. 3,
601-620.

\bibitem[I1]{Ishii-algebra}
{\sc Ishii, S.}
The arc space of a toric variety. {\em J. Algebra} {\bf  278}  (2004),  no. 2, 666--683.


 \bibitem[I2]{Ishii-crelle}
{\sc Ishii, S.}, Arcs, valuations and the Nash map, {\em J. reine
angew. Math.}, {\bf 588}, (2005), 71-92.

\bibitem[I3]{Ishii-fourier}
{\sc Ishii, S.}, The local Nash problem on arc families of
singularities, {\em Ann. Inst. Fourier (Grenoble)}. {\bf 56} (2006), no. 4, 1207--1224.

 \bibitem[J]{Jung}{\sc  Jung,  H.W.E.,}  Darstellung der Funktionen eines algebraischen K\"orpers zweier
 unabha\"angigen Ver\"anderlichen $x$, $y$ in der Umgebung einer stelle
  $x=a$, $y=b$, {\em  J.Reine Angew. Math.},  {\bf 133} (1908),
  289-314.




\bibitem[LJ-R]{LR}{\sc  Lejeune-Jalabert,  M.,  Reguera L\'opez, A. },
 Arcs and wedges on sandwiched surface singularities,
 {\em . Amer. J.
Math.} {\bf 121}, no. 6, (1999), 1191-1213.


\bibitem[L1]{Lipman2}{\sc  Lipman,  J.,} Introduction to resolution of singularities. {\em Proceedings of Symposia in Pure Mathematics}.
{\bf 29},  Amer. Math. Soc., Providence, R.I., 1975, 187--230.

\bibitem[L2]{Lipman}{\sc  Lipman,  J.,}
{\em Topological invariants of quasi-ordinary singularities,}
Memoirs of the American Mathematical Society  388, 1988.

\bibitem[Mo]{Morales} {\sc Morales M.}{The Nash problem on arcs for surface singularities}
arXiv:math/0609629.

\bibitem[N]{Nash} {\sc Nash, J. F. Jr.},  Arc structure of singularities. {\em A
celebration of John F. Nash, Jr. Duke Math. J. } {\bf 81} (1995),
no. 1, 31--38 (1996).

\bibitem[O]{Oda} {\sc Oda T.}; {\em Convex Bodies and algebraic
geometry}; Annals of Math. Studies(131), Springer-Verlag, 1988.


\bibitem[Pe]{Petrov}{\sc Petrov, P.}
Nash problem for stable toric varieties, arXiv: math.AG/0604432,
to appear in {\em Mathematische Nachrichten}.

\bibitem[Pl1]{Plenat1}{\sc Pl\'enat, C.}
\`A propos du probl\`eme des arcs de Nash. {\em Ann. Inst. Fourier (Grenoble)}
{\bf 55} (2005), no. 3, 805--823.


\bibitem[Pl2]{Plenat2} {\sc Pl\'enat, C.}
  R\'esolution du probl\`eme des arcs de Nash pour les points doubles rationnels $D\sb n$.
{\em C. R. Math. Acad. Sci. Paris} {\bf 340} (2005), no. 10, 747--750.


\bibitem[Pl3]{Plenat} {\sc Pl\'enat, C.} R\'esolution du probleme des arcs de Nash pour les
  points doubles rationnels.  Arxiv math.AG/0302188



\bibitem[Pl-PP1]{PPlenat}{\sc Plenat C., Popescu-Pampu, P.}
  A class of non-rational surface singularities with bijective Nash
  map.
  {\em Bull. Soc. Math. France} {\bf  134}  (2006),  no. 3, 383--394.



\bibitem[Pl-PP2]{PPlenat2}{\sc Plenat C., Popescu-Pampu, P.}
  Families of higher dimensional germs with bijective Nash map. Arxiv math.AG/0605566

\bibitem[PP]{PP}{\sc
Popescu-Pampu, P.} On the analytical invariance of the semigroups
of a
  quasi-ordinary hypersurface singularity.
{\em Duke Math. J.} {\bf 124}  (2004),  no. 1, 67-104.

\bibitem[Re1]{Reguera}{\sc Reguera-Lopez, A.J.} Families of arcs
on rational surface singularities. {\em Manuscr. Math.} {\bf 88},
1995, 321-333.

\bibitem[Re2]{Reguera2}{\sc Reguera-Lopez, A.J.}
Image of the Nash map in terms of wedges.
{\em C. R. Math. Acad. Sci. Paris} {\bf 338} (2004), no. 5, 385--390.


\end{thebibliography}
\end{document}